\documentclass[11pt]{amsart}
\usepackage{amssymb, latexsym}
\newtheorem{defin}{}
\newtheorem{saetze}[defin]{}
\newtheorem{conjec}[defin]{}
\newtheorem{lemmas}[defin]{}
\newtheorem{folger}[defin]{}
\newtheorem{bemerk}[defin]{}

\newenvironment{theorem}  {\begin{saetze}\it {\bf Theorem:}}{\end{saetze}}

\newcommand{\fillbox}{\mbox{$\bullet$}}

\begin{document}

\title[Theorem of Karrass and Solitar]
{A simple proof of a theorem of Karrass and Solitar}

\author[Delaram Kahrobaei]{Delaram Kahrobaei}
\thanks{The research of the author was
partially supported by the NSF, Grant 402350001.}
\author[Delaram Kahrobaei]{Delaram Kahrobaei}
\address{ Delaram Kahrobaei, Mathematical Institute, University of St
Andrews, North Haugh, St Andrews, Fife KY16 9SS Scotland, UK}
\email{delaram.kahrobaei@st-andrews.ac.uk}
\urladdr{http://www-groups.mcs.st-and.ac.uk/${\sim}$delaram/}

\begin{abstract} In this note we give a particularly short and simple proof
of the following theorem of Karrass and Solitar.  Let $H$ be a finitely
generated subgroup of a free group $F$  with infinite index $[F:H]$.
Then there is a nontrivial normal subgroup $N$ of $F$ such that
$N\cap H = \{1\}$.
\end{abstract}
\date{2004; to appear in AMS Proceedings, Contemporary Mathematics}

 \dedicatory{To Katalin Bencs\'{a}th} \maketitle
\section{Introduction}
As stated in the abstract, the purpose of this note is to present
a proof of a theorem of Karrass and Solitar \cite{ks} which is
particularly short and simple.  This result was also proved by
Ivanov and Schupp in \cite{is} who were apparently not aware of
\cite{ks}.
\newline
\noindent In \cite{M} and \cite{BMS02} (problem F22), P. Neumann
asks the following question: {\it Is it possible that the free
product $\{A \ast B; H = K\}$ with amalgamation, where $A$, $B$
are free groups of finite ranks, $H, K$ are finitely generated
subgroups of $A, B$, respectively, is a simple group?} In desire
of finding the answer to this problem, I arrived at the proof,
recorded below, to the following well-known theorem of Karrass and
Solitar. Note that as mentioned in \cite{is}, the theorem gives a
partial answer to the question of P. Neumann, that under the same
condition of the problem, $\{A \ast B; H=K\}$ is not simple
provided either of indices $[A:H], [B:K]$ is infinite. Needless to
mention that in \cite{BM97}, Burger and Mozes constructed an
example of an amalgamated product of finitely generated free
groups, with finitely generated amalgamated subgroup, which is
simple. The reader may want to recall the familiar sounding but
very different fact that if H is of {\it finite} index then $H
\cap K$ is non-trivial for {\it all} non-trivial $K < F$, (see pg.
20, \cite{LS77}).
\newline
\noindent In the proof of the theorem we use a theorem of M. Hall
\cite{mh} which states that a finitely generated subgroup of a
finitely generated free group is a free factor of a subgroup of
finite index.
\section{Main Theorem and Proof}
\begin{theorem} Let $H$ be a finitely generated subgroup of a free group $F$ with infinite index $[F:H]$.
Then there is a nontrivial normal subgroup $N$ of $F$ such that
$N\cap H = \{1\}$.
\end{theorem}
\begin{proof}
By a theorem of M. Hall \cite{mh}, $H$ is a free factor of a subgroup
$K= H*Q$ which has finite index  in $F$.
Since $H$ has infinite index, $Q$ must be nontrivial.
There are only finitely many conjugates of $K$ in $F$
and their intersection $C$ is a normal
subgroup of $F$ of finite index with $C\leq K =H*Q$.

By the Kurosh subgroup theorem $C= (C\cap H) * J$ for some subgroup
$J$.  Since $C$ has finite index in $K$,  $J$ is nontrivial.  Let $L=\mbox{gp}_C(J)$
be the normal closure of $J$ in $C$.  Then $L\cap H= L\cap(C\cap H) = \{1\}$.

Let $b_1=1,b_2,\ldots,b_n$ be a complete system of coset representatives for
$C$ in $F$.  Conjugation by $b_i$ induces an automorphism of $C$ and so each
$L^{b_i}$ is normal in $C$.  The intersection $N = \cap_{i=1}^n L^{b_i}$
is clearly normal in $F$.  Also $N$ contains the nontrivial subgroup
$[L^{b_1},L^{b_2},\ldots,L^{b_n}]$ because each $L^{b_i}$ is normal in $C$.
Finally $N\cap H =\{1\}$ since $N \leq L$,
so $N$ is the desired nontrivial normal subgroup.
\end{proof}
\subsection*{Acknowledgement}
I would like to thank G.Baumslag, C.F.Miller and the anonymous
referee for their suggestions.


\begin{thebibliography}{999}

\bibitem{M}
{\em New York Group Theory Cooperative at CCNY},
http://grouptheory.org.

\bibitem{BMS02} G.~Baumslag and A.~G.~Myasnikov and V.~Shpilrain,
{\em Open problems in combinatorial group theory. Second edition.
Combinatorial and geometric group theory}, Contemp. Math. {\bf
296} (2002), 1--38.

\bibitem{mh} M.~Hall~Jr.,
{\em Subgroups of finite index in free groups},
Canad. J. Math. {\bf 1} (1949), 187--190.

\bibitem{is} S.~V.~Ivanov and P.~E.~Schupp,
{\em A remark on finitely generated subgroups of a free group}, in
``Algorithmic problem in groups and semigroups (Lincoln NE,
1998)", (J.-C. Birget, S. Margolis, J. Meakin, M. Sapir, Eds.),
Trends in Mathematics, Birkhauser, Boston (2000), 139-142.

\bibitem{ks} A.~Karrass and D.~Solitar,
{\em On finitely generated subgroups of a free group},
Proc. Amer. Math. Soc. {\bf 22} (1969), 209-213.

\bibitem{LS77} R.~C.~Lyndon and P.~E.~Schupp,
{\em Combinatorial Group Theory}, in "Classics in Mathematics",
Springer-Verlag, Berlin (1977).

\bibitem{BM97} M.~Burger and S.~Mozes,
{\em Finitely presented simple groups and products of trees}, C.
R. Acad. Sci. Paris Serie I Math. {\bf 7} (1997), 747-752.


\end{thebibliography}
\end{document}